\documentclass[a4paper,11pt]{article}

\usepackage{pdfpages}
\usepackage{graphicx}
\usepackage{mathtools}
\usepackage{verbatim}
\usepackage{amssymb}
\usepackage{bm}
\usepackage{enumitem}
\usepackage{cases}
\usepackage{latexsym}
\usepackage{amsthm}
\usepackage{amsfonts}
\usepackage{bbm}
\usepackage{etoolbox}
\usepackage[top=25mm, bottom=25mm, left=28mm, right=28mm]{geometry}
\usepackage{cite}
\usepackage{changepage}
\usepackage{fancybox}
\usepackage{tikz}

\newtheorem{theorem}{Theorem}[section]
\AfterEndEnvironment{theorem}{\noindent\ignorespaces}
\newtheorem{corollary}[theorem]{Corollary}
\AfterEndEnvironment{corollary}{\noindent\ignorespaces}
\newtheorem{lemma}[theorem]{Lemma}
\AfterEndEnvironment{lemma}{\noindent\ignorespaces}

\AfterEndEnvironment{proposition}{\noindent\ignorespaces}

\AfterEndEnvironment{question}{\noindent\ignorespaces}
\newtheorem{conjecture}[theorem]{Conjecture}
\AfterEndEnvironment{conjecture}{\noindent\ignorespaces}

\AfterEndEnvironment{problem}{\noindent\ignorespaces}

\theoremstyle{definition}
\newtheorem{definition}[theorem]{Definition}
\AfterEndEnvironment{definition}{\noindent\ignorespaces}

\theoremstyle{remark}

\AfterEndEnvironment{example}{\noindent\ignorespaces}

\AfterEndEnvironment{remark}{\noindent\ignorespaces}

\AfterEndEnvironment{notation}{\noindent\ignorespaces}

\newenvironment{proof2}{\proof[\textnormal{\textbf{Proof.}}]}{\qed}

\usepackage{chngcntr}
\usepackage{apptools}
\AtAppendix{\counterwithin{theorem}{section}}

\def\e{\varepsilon}

\def\N{\mathbb{N}}

\def\d{\delta}

\def\s{\subset}
\def\ex{\mathrm{ex}}

\def\cP{\mathcal{P}}
\def\cF{\mathcal{F}}

\begin{document}

\title{The extremal number of longer subdivisions}

\author{Oliver Janzer\thanks{Department of Pure Mathematics and Mathematical Statistics, University of Cambridge, United Kingdom.
E-mail: {\tt oj224@cam.ac.uk}.}}

\date{}

\maketitle

\begin{abstract}
	
For a multigraph $F$, the $k$-subdivision of $F$ is the graph obtained by replacing the edges of $F$ with pairwise internally vertex-disjoint paths of length $k+1$. Conlon and Lee conjectured that if $k$ is even, then the $(k-1)$-subdivision of any multigraph has extremal number $O(n^{1+\frac{1}{k}})$, and moreover, that for any simple graph $F$ there exists $\e>0$ such that the $(k-1)$-subdivision of $F$ has extremal number $O(n^{1+\frac{1}{k}-\e})$. In this paper, we prove both conjectures.

\end{abstract}

\section{Introduction}

For a multigraph $F$, a \emph{subdivision} of $F$ is a graph obtained by replacing the edges of $F$ with pairwise internally vertex-disjoint paths of arbitrary lengths. The $k$-subdivision of $F$ is the graph obtained by replacing the edges of $F$ with pairwise internally vertex-disjoint paths of length $k+1$, and is denoted by $F^k$.

Many researchers have studied the problem of estimating the number of edges needed in a graph $G$ on $n$ vertices to guarantee that it contains as a subgraph a subdivided copy of a fixed graph. The first result in this direction is due to Mader \cite{Ma67} who proved that for any graph $F$ there exists a constant $c_F=c$ such that if an $n$-vertex graph $G$ contains at least $cn$ edges, then $G$ contains a subdivision of $F$ as a subgraph. In this result the size of the subdivided graph can grow with $n$, which is necessary since an $n$-vertex graph with $cn$ edges need not contain a cycle of bounded length.

Answering a question of Erd\H os about planar subgraphs \cite{Er71}, Kostochka and Pyber \cite{KP88} proved that any $n$-vertex graph with at least $4^{t^2}n^{1+\e}$ edges contains a subdivided $K_t$ with at most $\frac{7t^2\log t}{\e}$ vertices. This is the first result that guarantees a subdivided $K_t$ of bounded size.

For a family $\cF$ of graphs, we let $\ex(n,\cF)$ be the maximum number of edges in an $n$-vertex graph not containing any $F\in \cF$ as a subgraph. When $\cF=\{F\}$, we write $\ex(n,F)$ for the same function.

Let $\cF_{t,k}$ be the family of graphs that can be obtained by replacing the edges of $K_t$ with pairwise internally vertex-disjoint paths of length at most $k$. Jiang \cite{J11} proved that for any $t\in \N$ and any $0<\e<1/2$, we have $\ex(n,\cF_{t,\lceil 10/\e\rceil})=O(n^{1+\e})$. Here the asymptotic notation means that $n\rightarrow \infty$ and other parameters are constant. We follow the same convention throughout the paper.

Note that Jiang's result improves that of Kostochka and Pyber in two ways. Firstly, any $F\in \cF_{t,\lceil 10/\e\rceil}$ has at most $\frac{ct^2}{\e}$ vertices, so a $\log$ factor is saved. Secondly, the edges in Jiang's theorem are replaced by uniformly short paths not depending on $t$. However, they can still have different lengths. The next result of Jiang and Seiver guarantees a subdivided $K_t$ with prescribed path lengths.

\begin{theorem}[Jiang--Seiver \cite{JS12}] \label{thmJS}
	For any $t\in \N$ and any even $k\in \N$, $$\ex(n,K_t^{k-1})=O(n^{1+\frac{16}{k}}).$$
\end{theorem}

Note that if $k$ is odd, then $K_t^{k-1}$ is not a bipartite graph, so $\ex(n,K_t^{k-1})=\Theta(n^2)$.

Conlon and Lee conjectured that the following two strengthenings hold.

\begin{conjecture}[Conlon--Lee \cite{CL18}] \label{multigraphconj}
	Let $F$ be a multigraph and let $k\geq 2$ be even. Then $$\ex(n,F^{k-1})=O(n^{1+\frac{1}{k}}).$$
\end{conjecture}

\begin{conjecture}[Conlon--Lee \cite{CL18}] \label{simplegraphconj}
	Let $F$ be a simple graph and let $k\geq 2$ be even. Then there exists some $\e>0$ such that $$\ex(n,F^{k-1})=O(n^{1+\frac{1}{k}-\e}).$$
\end{conjecture}

In the case $k=2$, Conjecture \ref{multigraphconj} follows from the $r=2$ case of a result of F\"uredi \cite{Fu91} and Alon, Krivelevich and Sudakov \cite{AKS03}, which states that any bipartite graph with maximum degree at most $r$ on one side has extremal number $O(n^{2-1/r})$. The $k=2$ case of Conjecture~\ref{simplegraphconj} was proved by Conlon and Lee \cite{CL18}, and improved bounds were given by the author \cite{Ja18}.

Very recently, Conlon, Janzer and Lee proved Conjecture \ref{simplegraphconj} for every bipartite graph~$F$.

\begin{theorem}[Conlon--Janzer--Lee \cite{CJL19}] \label{thmCJL}
	Let $F$ be a simple bipartite graph and let $k\geq 1$. Then there exists some $\e>0$ such that $$\ex(n,F^{k-1})=O(n^{1+\frac{1}{k}-\e}).$$
\end{theorem}

As a simple corollary, they significantly improved the bound in Theorem \ref{thmJS}.

\begin{theorem}[Conlon--Janzer--Lee \cite{CJL19}]
	Let $F$ be a simple graph and let $k\geq 2$ be even. Then there exists some $\e>0$ such that $$\ex(n,F^{k-1})=O(n^{1+\frac{2}{k}-\e}).$$
\end{theorem}

In this paper, we prove both Conjecture \ref{multigraphconj} and Conjecture \ref{simplegraphconj}.

\begin{theorem} \label{multisub}
	Let $F$ be a multigraph and let $k\geq 2$ be even. Then $$\ex(n,F^{k-1})=O(n^{1+\frac{1}{k}}).$$
\end{theorem}

\begin{theorem} \label{simplesub}
	Let $F$ be a simple graph and let $k\geq 2$ be even. Then there exists some $\e>0$ such that $$\ex(n,F^{k-1})=O(n^{1+\frac{1}{k}-\e}).$$
\end{theorem}

Note that these results are tight. Indeed, by a result of Conlon \cite{C18}, the Theta graph $\theta_{k,\ell}$ has extremal number $\Theta(n^{1+1/k})$ for all $\ell\geq \ell_0(k)$, showing that Theorem \ref{multisub} is tight. Moreover, Erd\H os-R\'enyi random graphs show that $\ex(n,K_t^{k-1})=\Omega(n^{1+1/k-c_{k,t}})$ where $c_{k,t}\rightarrow 0$ as $t\rightarrow \infty$, so Theorem \ref{simplesub} is also tight.

The rest of the paper is organised as follows. In Section \ref{secoverall}, we introduce some of the key definitions and give the high-level structure of the proof, with the key technical lemmas deferred to Sections \ref{secshortpaths} and \ref{seclongpaths}.

\section{The high-level structure of the proof} \label{secoverall}

A graph $G$ is called $K$-almost-regular if $\max_{v\in V(G)} d(v)\leq K\min_{v\in V(G)} d(v)$, where $d(v)$ is the degree of vertex $v$. The following lemma, which is a small modification of a result proved by Erd\H os and Simonovits \cite{ES70}, allows us to restrict our attention to almost regular host graphs.

\begin{lemma}[Jiang--Seiver \cite{JS12}] \label{lemmaJS}
	Let $\e,c$ be positive reals, where $\e<1$ and $c\geq 1$. Let $n$ be a positive integer that is sufficiently large as a function of $\e$. Let $G$ be a graph on $n$ vertices with $e(G)\geq cn^{1+\e}$. Then $G$ contains a $K$-almost-regular subgraph $G'$ on $m\geq n^{\frac{\e-\e^2}{2+2\e}}$ vertices such that $e(G')\geq \frac{2c}{5}m^{1+\e}$ and $K=20\cdot 2^{\frac{1}{\e^2}+1}$.
\end{lemma}

Using this lemma, Theorem \ref{multisub} and Theorem \ref{simplesub} reduce to the following two statements, respectively. For notational convenience, we have dropped the assumption that $k$ is even, and replaced $k$ by $2k$.

\begin{theorem} \label{multisubreduced}
	Let $F$ be a multigraph and let $k\geq 1$. Suppose that $G$ is a $K$-almost-regular graph on $n$ vertices with minimum degree $\d=\omega(n^{\frac{1}{2k}})$. Then, for $n$ sufficiently large, $G$ contains a copy of $F^{2k-1}$.
\end{theorem}

\begin{theorem} \label{simplesubreduced}
	Let $F$ be a simple graph and let $k\geq 1$. Then there exists $\e>0$ with the following property. Suppose that $G$ is a $K$-almost-regular graph on $n$ vertices with minimum degree $\d=\omega(n^{\frac{1}{2k}-\e})$. Then, for $n$ sufficiently large, $G$ contains a copy of $F^{2k-1}$.
\end{theorem}

From now on we let $F$ be an arbitrary fixed multigraph and write $H=F^{2k-1}$. Moreover, throughout the paper we tacitly assume that $n$ is sufficiently large.

The next definition was introduced in \cite{CJL19}, and was used to prove Theorem \ref{thmCJL}.

\begin{definition}
	Let $L$ be a positive real and let $f(\ell,L)=L^{5^\ell}$ for $1\leq \ell\leq 2k$.
	We recursively define the notions of \emph{$L$-admissible} and \emph{$L$-good paths of length $\ell$} in a graph.
	Any path of length $1$ is both $L$-admissible and $L$-good. For $2\leq \ell\leq 2k$, we say a path $P=v_0v_1\dots v_\ell$ is $L$-admissible if every proper subpath of $P$ is $L$-good, i.e., $v_iv_{i+1}\dots v_j$ is $L$-good for every $(i,j)\neq (0,\ell)$. The path $P$ is $L$-good if it is $L$-admissible and the number of $L$-admissible paths of length $\ell$ between $v_0$ and $v_\ell$ is at most $f(\ell,L)$.
\end{definition}

The next lemma will be used several times later.

\begin{lemma} \label{manydisjoint}
	Let $\ell\geq 2$ and let $L>\ell$. If a path $P=v_0\dots v_\ell$ is $L$-admissible, but not $L$-good, then there exist at least $L$ pairwise internally vertex-disjoint paths of length $\ell$ from $v_0$ to $v_\ell$.
\end{lemma}

\begin{proof2}
	Take a maximal set of pairwise internally vertex-disjoint paths of length $\ell$ from $v_0$ to $v_\ell$ and assume that it consists of fewer than $L$ paths. These paths contain at most $L(\ell-1)$ internal vertices in total and any path of length $\ell$ between $v_0$ and $v_\ell$ intersects at least one of these vertices. Since there are at least $L^{5^\ell}$ $L$-admissible paths of length $\ell$ between $v_0$ and $v_\ell$, it follows by pigeon hole that there exist some $1\leq i\leq \ell-1$ and some $x\in V(G)$ such that there are at least $\frac{L^{5^{\ell}}}{(\ell-1)L(\ell-1)}$ $L$-admissible paths of the form $u_0u_1\dots u_\ell$ with $u_0=v_0,u_i=x,u_\ell=v_\ell$. Observe that $\frac{L^{5^{\ell}}}{(\ell-1)L(\ell-1)}>L^{5^i}L^{5^{\ell-i}}$, so either there are more than $L^{5^{i}}$ $L$-good paths of length $i$ between $v_0$ and $x$ or there are more than $L^{5^{\ell-i}}$ $L$-good paths of length $\ell-i$ between $x$ and $v_\ell$. In either case, we contradict the definition of an $L$-good path.
\end{proof2}


\medskip

Our strategy will be to prove that, roughly speaking, in any almost regular $H$-free graph there are many good paths of length $2k$. As we will see in Section \ref{secshortpaths}, the techniques in \cite{CJL19} can be easily applied to prove this for paths of length $k$. The novelty of this paper is the machinery that allows us to extend this to longer paths, using very different techniques. This is given in Section \ref{seclongpaths}, where we prove the following lemma.

\begin{lemma} \label{good2k}
	Let $G$ be an $H$-free $K$-almost-regular graph on $n$ vertices with minimum degree~$\d\geq L^{100^k|V(H)|}$, and let $S\s V(G)$. Then, provided that $L$ is sufficiently large compared to $|V(H)|$ and $K$, $|S|=\omega(\frac{n}{\d^{1/2}})$ and $|S|=\omega(\frac{n}{L^{1/2}})$, the number of $L$-good paths of length $2k$ with both endpoints in $S$ is $\Omega(\frac{|S|^2\d^{2k}}{n})$.
\end{lemma}

Note that in this result and everywhere else in the paper, the asymptotic notation $\Omega$ allows the implied constant to depend on $k,|V(H)|$ and $K$, which are thought of as constants, while $\d$ and $L$ are functions of $n$.

With Lemma \ref{good2k} in hand, the proof of Theorem \ref{multisubreduced} is immediate.

\begin{proof}[\textnormal{\textbf{Proof of Theorem \ref{multisubreduced}}}]
	Suppose that $G$ does not contain $H=F^{2k-1}$ as a subgraph. Since $\d=\omega(n^{\frac{1}{2k}})$, we may choose $L$ with $L=\omega(1)$, $L^{100^k|V(H)|}\leq \d$ and $n^2f(2k,L)=o(n\d^{2k})$. Then we may apply Lemma \ref{good2k} with $S=V(G)$ to get that the number of $L$-good paths of length $2k$ in $G$ is $\Omega(n\d^{2k})$, which is $\omega(n^2f(2k,L))$. However, by the definition of $L$-goodness, between any two vertices there can be at most $f(2k,L)$ such paths, which is a contradiction.
\end{proof}

The proof of Theorem \ref{simplesubreduced} is slightly more complicated, and it uses ideas from \cite{Ja18}.

\begin{proof}[\textnormal{\textbf{Proof of Theorem \ref{simplesubreduced}}}]
	Firstly note that $F$ is a subgraph of $K_t$ for some $t$, so it suffices to prove the result for $F=K_t$. Let $\e>0$ be sufficiently small, to be specified, and let $G$ be a $K$-almost-regular graph on $n$ vertices with minimum degree $\d=\omega(n^{\frac{1}{2k}-\e})$. Assume that $G$ does not contain a copy of $H=F^{2k-1}$.
	
	For vertices $u,v\in V(G)$, let us write $u\sim v$ if there is a path of length $2k$ between $u$ and $v$. Also, let us say that $u$ and $v$ are \emph{distant} if for every $1\leq i\leq 4k-2$, the number of walks of length $i$ between $u$ and $v$ is at most $\d^{i-2k+1/2}$. Observe that for any $u\in V(G)$ the number of walks of length $i$ starting from $u$ is at most $(K\d)^{i}$, so the number of vertices $v\in V(G)$ for which there are at least $\d^{i-2k+1/2}$ walks of length $i$ from $u$ to $v$ is at most $\frac{(K\d)^{i}}{\d^{i-2k+1/2}}=K^{i}\d^{2k-1/2}$. Thus, the number of $v\in V(G)$ for which $u$ and $v$ are not distant is $O(\d^{2k-1/2})$.
	
	Define $c_0=\e$ and $c_{\ell+1}=(3\cdot 5^{2k}+1)c_\ell+2k\e$ for $0\leq \ell\leq t-1$. Assume that $\e$ is small enough so that
	\begin{equation}
		3\cdot 100^k|V(H)|\cdot c_\ell\leq \frac{1}{2k}-\e \label{eqn100tok}
	\end{equation}
	for all $0\leq \ell\leq t$.
	Then in particular $c_\ell\leq \frac{1}{4k}-\e/2$ holds for all $0\leq \ell\leq t$. For future reference, note that then \begin{equation}
	n^{c_\ell}\leq n^{\frac{1}{4k}-\e/2}= o(\delta^{1/2}). \label{deltaandc}
	\end{equation}
	
	\medskip
	
	\noindent \emph{Claim.} For any $0\leq \ell\leq t$, there exist distinct vertices $x_1,\dots,x_\ell\in V(G)$ and a set $S_\ell\s V(G)$ such that
	\begin{enumerate}[label=(\roman*)]
		\item there is a copy of $K_\ell^{2k-1}$ in $G$ with the vertices of the subdivided $K_\ell$ being $x_1,\dots,x_\ell$
		\item $x_i\sim y$ for every $1\leq i\leq \ell$ and every $y\in S_\ell$
		\item $|S_\ell|=\Omega(n^{1-c_\ell})$ and
		\item $x_i$ and $x_j$ are distant for every $1\leq i<j\leq \ell$.
	\end{enumerate}
	
	\medskip
	
	Note that in particular for $\ell=t$, condition (i) guarantees the existence of a subgraph $K_t^{2k-1}$, so it suffices to prove the claim.
	
	\medskip
	
	\noindent \emph{Proof of Claim.} We proceed by induction on $\ell$. For $\ell=0$, we may take $S_0=V(G)$. Assume now that we have verified the claim for $\ell$.
	
	Suppose that for some $y\in S_\ell$ there exist $1\leq i<j\leq \ell$ and two paths of length $2k$, one (called $P_i$) from $x_i$ to $y$ and one (called $P_j$) from $x_j$ to $y$, which share a vertex other than $y$. Let they intersect at some vertex $z\neq y$. Now let the subpath of $P_i$ between $x_i$ and $z$ have length $\alpha$ and let the subpath of $P_j$ between $x_j$ and $z$ have length $\beta$. Then there is a walk of length $\alpha+\beta$ from $x_i$ to $x_j$ through $z$. Moreover, there is a path of length $2k-\alpha$ from $z$ to $y$. Observe that $2k-\alpha\leq 4k-(\alpha+\beta)-1$.
	
	Let $Y$ be the set of $y\in S_\ell$ for which there exist some $1\leq i<j\leq \ell$ and a walk $W$ of length $\gamma \leq 4k-2$ between $x_i$ and $x_j$ such that for some vertex $w$ on $W$ the distance of $y$ from $w$ is at most $4k-\gamma-1$.
	By condition (iv), there are at most $\d^{\gamma-2k+1/2}$ walks of length $\gamma$ between any $x_i$ and $x_j$ so there are $O(\d^{\gamma-2k+1/2})$ vertices appearing in at least one of these walks. Therefore the number of vertices at distance at most $4k-\gamma-1$ from at least one of these vertices is $O(\d^{\gamma-2k+1/2}\cdot \d^{4k-\gamma-1})=O(\d^{2k-1/2})$. That is, $|Y|=O(\d^{2k-1/2})$.
	
	Notice that by the discussion above, for any $y\in S_\ell\setminus Y$ and any $i\neq j$, a path of length $2k$ from $x_i$ to $y$, and a path of length $2k$ from $x_j$ to $y$ have no common vertex other than $y$. Thus, by condition (ii) there exist $\ell$ paths of length $2k$, one from each $x_i$ to $y$ which are pairwise vertex-disjoint apart from at $y$. Moreover, these paths are also vertex-disjoint from the paths forming the $K_\ell^{2k-1}$ guaranteed by condition (i), apart from the trivial intersections at $x_1,\dots,x_\ell$ (else, there is a path of length at most $2k-1$ from $y$ to a point on a path of length $2k$ between some $x_i$ and $x_j$, which contradicts the fact that $y\not \in Y$). Thus, for any $y\in S_\ell\setminus Y$ there is a copy of $K_{\ell+1}^{2k-1}$ in $G$ with the vertices of the subdivided $K_{\ell+1}$ being $x_1,\dots,x_\ell,y$.
	
	Let $Z$ be the set of $z\in S_\ell$ which are not distant to $x_i$ for at least one $1\leq i\leq \ell$. By the second paragraph in this proof, $|Z|=O(\d^{2k-1/2})$.
	
	Let $S'_\ell=S_\ell\setminus (Y\cup Z)$. Recall that $|Y|=O(\d^{2k-1/2})$. Note that if $\d=\omega(n^{\frac{1}{2k}})$, then, by Theorem \ref{multisubreduced}, $G$ contains $H$ as a subgraph, so we may assume that $\d=O(n^{\frac{1}{2k}})$. Then $\d^{2k-1/2}=O(\frac{n}{\d^{1/2}})$, which is $o(n^{1-c_\ell})$ by equation (\ref{deltaandc}). Thus, $|Y\cup Z|=o(n^{1-c_\ell})$ and so $|S'_\ell|= \Omega(n^{1-c_\ell})$.
	
	Let $L=n^{3c_\ell}$. Then, by equation (\ref{eqn100tok}), we have $L^{100^k|V(H)|}\leq n^{\frac{1}{2k}-\e}=o(\d)$.  Moreover, by equation (\ref{deltaandc}), we have $n^{1-c_\ell}=\omega(\frac{n}{\d^{1/2}})$, and by the definition of $L$, we have $n^{1-c_\ell}=\omega(\frac{n}{L^{1/2}})$. Hence, by Lemma \ref{good2k}, the number of $L$-good paths of length $2k$ with both endpoints in $S'_\ell$ is $\Omega(\frac{|S'_\ell|^2\d^{2k}}{n})$. Between any two vertices in $S'_\ell$ there are at most $f(2k,L)$ $L$-good paths of length $2k$, so the number of pairs $(z,y)\in S'_\ell\times S'_\ell$ with $z\sim y$ is $\Omega(\frac{|S'_\ell|^2\d^{2k}}{nf(2k,L)})$. Thus, there exists some $x_{\ell+1}\in S'_\ell$ such that the number of $y\in S'_{\ell}$ with $x_{\ell+1}\sim y$ is $\Omega(\frac{|S'_\ell|\d^{2k}}{nf(2k,L)})\geq\Omega(\frac{n^{1-c_\ell-2k\e}}{L^{5^{2k}}})=\Omega(n^{1-c_\ell-2k\e-3c_\ell 5^{2k}})=\Omega(n^{1-c_{\ell+1}})$. Set $S_{\ell+1}$ to be the set of these $y\in S'_{\ell}$, and note that properties (i)-(iv) are satisfied for $\ell+1$.
\end{proof}


\section{Short paths} \label{secshortpaths}

Our aim in this section is to prove the following lemma.

\begin{lemma} \label{kpathsgood}
	Let $G$ be an $H$-free $K$-almost-regular graph on $n$ vertices with minimum degree~$\d\geq L^{100^k|V(H)|}$. Then, provided that $L$ is sufficiently large compared to $|V(H)|$ and $K$, the number of paths of length $k$ that are not good is $O(\frac{n\d^k}{L})$.
\end{lemma}

The proof of this is almost identical to that of Lemma 6.4 in \cite{CJL19}, nevertheless we include it here for completeness and since some minor details need to be modified.

The next definition is for notational convenience.

\begin{definition}
	A pair of distinct vertices $(x,y)$ in $G$ is said to be $(\ell,L)$-\emph{bad} for some $2 \leq \ell \leq 2k$ and some $L$ if there is an $L$-admissible, but not $L$-good, path of length $\ell$ from $x$ to $y$.
\end{definition}

In what follows, for $v\in V(G)$, we shall write $\Gamma_i(v)$ for the set of vertices $u\in V(G)$ for which there exists a path of length $i$ from $v$ to $u$ and write $N(v)=\Gamma_1(v)$. The next lemma is a slight variant of Lemma 6.7 from \cite{CJL19}. 

\begin{lemma} \label{findstructure}
	Let $2\leq \ell\leq k$ and $1\leq i\leq \ell$. Let $G$ be a $K$-almost-regular graph on $n$ vertices with minimum degree $\d>0$. Let $X,Y,Z\s V(G)$ be such that $|Z|\leq L^{1/10},|Y|\leq (K\d)^{\ell-1}$ and, for any $x\in X$, the number of $y\in Y$ such that $(x,y)$ is $(\ell,L)$-bad is as at least $\frac{(K\d)^{\ell-1}}{f(\ell-1,L)^2}$. Then, provided that $L$ is sufficiently large compared to $k$ and $K$, there exist a path of length $2i$ in $G$, disjoint from $Z$, whose endpoints form a set $R\s Y$, and a subset $X'\s X$ such that $|X'|\geq |X\setminus Z|/(16f(\ell-1,L)^2)$ and $(x',r)$ is $(\ell,L)$-bad for every $x'\in X'$ and $r\in R$.
\end{lemma}

\begin{proof2}
	After replacing $X$ by $X\setminus Z$, we may assume $X\cap Z=\emptyset$. Let $Y'$ be the set of those $y\in Y$ for which the number of $x\in X$ such that $(x,y)$ is $(\ell,L)$-bad is at least $\frac{|X|}{2f(\ell-1,L)^2}$. 
	Then the number of $(x,y)\in X\times (Y\setminus Y')$ which are $(\ell,L)$-bad is at most $\frac{|X||Y|}{2f(\ell-1,L)^2}\leq \frac{|X|(K\d)^{\ell-1}}{2f(\ell-1,L)^2}$, so the number of $(x,y)\in X\times Y'$ which are $(\ell,L)$-bad is at least $\frac{|X|(K\d)^{\ell-1}}{2f(\ell-1,L)^2}$. Now there exists some $x^*\in X$ such that there are at least $\frac{(K\d)^{\ell-1}}{2f(\ell-1,L)^2}$ choices $y\in Y'$ for which $(x^*,y)$ is $(\ell,L)$-bad. If a pair $(x^*,y)$ is $(\ell,L)$-bad, then there are at least $f(\ell,L)$ paths of length $\ell$ from $x^*$ to $y$. Hence, there are at least $\frac{(K\d)^{\ell-1}}{2f(\ell-1,L)^2}\cdot f(\ell,L)= \Omega(f(\ell-1,L)^3\d^{\ell-1})$ paths of length $\ell$ starting at $x^*$ and ending in $Y'$. 
    
    The number of such paths intersecting $Z$ is at most $|Z|\ell(K\d)^{\ell-1}$. Indeed, there are at most $|Z|$ choices for the element of $Z$ in the path, at most $\ell$ choices for its position in the path and, given a fixed choice for these, at most $(K\d)^{\ell-1}$ choices for the other $\ell-1$ vertices in the path. (Note that as $X\cap Z=\emptyset$, the vertex in $Z$ is not $x^*$.) But $|Z|\ell(K\d)^{\ell-1}\leq L^{1/10}\ell K^{\ell-1}\d^{\ell-1}$, so, for $L$ sufficiently large there are $\Omega(f(\ell-1,L)^3\d^{\ell-1})$ paths of length $\ell$ starting at $x^*$ and ending in $Y'$ that avoid $Z$. Moreover, since $|\Gamma_{\ell-i}(x^*)|\leq (K\d)^{\ell-i}$, it follows that there exists some $u\in \Gamma_{\ell-i}(x^*)$ such that there are $\Omega(f(\ell-1,L)^3\d^{i-1})$ paths of length $i$ from $u$ to $Y'$, all avoiding $Z$.
    
    Take now a maximal set of such paths which are pairwise vertex-disjoint apart from at~$u$. We claim that there are $\Omega(f(\ell-1,L)^3)$ such paths. Suppose otherwise. Then all the $\Omega(f(\ell-1,L)^3\d^{i-1})$ paths of length $i$ from $u$ to $Y'$ intersect a certain set of size $o(f(\ell-1,L)^3)$ not containing $u$. But there are $o(f(\ell-1,L)^3)\d^{i-1}$ such paths, which is a contradiction.
    
    So we have $r=\Omega(f(\ell-1,L)^3)$ paths $P_1,\dots,P_r$ of length $i$ from $u$ to $Y'$ which are pairwise vertex-disjoint except at $u$ and avoid $Z$. Let the endpoints of these paths be $y_1,\dots,y_r$. Since $y_j\in Y'$ for all $j$, the number of pairs $(x,y_j)$ with $x\in X$ which are $(\ell,L)$-bad is at least $\frac{r|X|}{2f(\ell-1,L)^2}$. Therefore, by Jensen's inequality, for an average $x\in X$ there are at least ${r/(2f(\ell-1,L)^2) \choose 2}$ choices $1\leq j_1<j_2\leq r$ such that both $(x,y_{j_1})$ and $(x,y_{j_2})$ are $(\ell,L)$-bad. Since ${r/(2f(\ell-1,L)^2) \choose 2}\geq (\frac{1}{4f(\ell-1,L)^2})^2{r \choose 2}$, there exist $1\leq j_1<j_2\leq r$ such that the set $$X'=\{x\in X: (x,y_{j_1}) \text{ and } (x,y_{j_2}) \text{ are } (\ell,L)\text{-bad}\}$$ has size at least $|X|/(4f(\ell-1,L)^2)^2$. We can now take $R=\{y_{j_1},y_{j_2}\}$, and the union of the paths $P_{j_1}$ and $P_{j_2}$ is a suitable path of length $2i$.
\end{proof2}

\medskip

The following lemma is a small modification of Lemma 6.8 from \cite{CJL19}. 

\begin{lemma}
	Let $G$ be an $H$-free $K$-almost-regular graph on $n$ vertices with minimum degree~$\d\geq L^{100^k|V(H)|}$. Let $2\leq \ell\leq k$ and any $v\in V(G)$. Then, provided that $L$ is sufficiently large compared to $|V(H)|$ and $K$, the number of $L$-admissible, but not $L$-good, paths of the form $v_0vv_2v_3\dots v_\ell$ is at most $\frac{2(K\d)^\ell}{f(\ell-1,L)}$.
\end{lemma}

\begin{proof2}
	
	
	Suppose otherwise. Let $Y=\Gamma_{\ell-1}(v)$ and note that $|Y|\leq (K\d)^{\ell-1}$. For any $x\in N(v)$ and any $y\in Y$, the number of $L$-admissible paths of the form $xvv_2\dots v_{\ell-1}y$ is at most $f(\ell-1,L)$. Indeed, in any such path, the subpath $vv_2v_3\dots v_{\ell-1}y$ is $L$-good, and for any fixed $y\in Y$ there are at most $f(\ell-1,L)$ such $L$-good paths. Hence, by assumption, the number of pairs $(x,y)\in N(v)\times Y$ such that there is an $L$-admissible, but not $L$-good, path of the form $xvv_1\dots v_{\ell-1}y$ is at least $\frac{2(K\d)^\ell}{f(\ell-1,L)^2}\geq \frac{2|N(v)|(K\d)^{l-1}}{f(\ell-1,L)^2}$. By definition, any such pair $(x,y)$ is $(\ell,L)$-bad. Let $X$ consist of those $x\in N(v)$ for which there are at least $\frac{(K\d)^{\ell-1}}{f(\ell-1,L)^2}$ choices of $y\in Y$ such that $(x,y)$ is $(\ell,L)$-bad. Then the number of pairs $(x,y)\in X\times Y$ which are $(\ell,L)$-bad is at least $\frac{|N(v)|(K\d)^{\ell-1}}{f(\ell-1,L)^2}$, and so $|X|\geq \frac{|N(v)|}{f(\ell-1,L)^2}\geq \frac{\d}{f(\ell-1,L)^2}$.
	
	Our aim now is to find a copy of $H$ in $G$, which will yield a contradiction. Write $k=j\ell+i$ with $1\leq i\leq \ell$.
	
	Note that if $L$ is sufficiently large, then $$|X|\geq \frac{\d}{f(\ell-1,L)^2}\geq \frac{L^{100^k|V(H)|}}{f(\ell-1,L)^2}\geq \frac{f(\ell-1,L)^{20|V(H)|}}{f(\ell-1,L)^2} \geq  2L(16f(\ell-1,L)^2)^{2|V(H)|},$$
	so we may apply Lemma \ref{findstructure} repeatedly $|E(F)|+|V(H)|\leq 2|V(H)|$ times and still get a set $X'$ of size at least $L$. Thus, we find disjoint paths $P_e$ of length $2i$ for every $e\in E(F)$ whose endpoint sets are $R_e\s Y$, and sets $X_{\text{final}}\s X$ and $U\s Y$ with $|X_{\text{final}}|=|U|=|V(H)|$ such that $V(P_e), X_{\text{final}}$ and $U$ are pairwise disjoint and any pair $(x,y)$ with $x\in X_{\text{final}}$ and $y\in U\cup \bigcup_{e\in E(F)} R_e$ is $(\ell,L)$-bad.
	For $e\in E(F)$, let $y_{e_{-k}}y_{e_{-k+1}}\dots y_{e_k}$ be the path of length~$2k$ replacing the edge $e$.
	
	A copy of $H$ in $G$ can now be constructed as follows. For each $e\in E(F)$, map the path $y_{e_{-i}}y_{e_{-i+1}}\dots y_{e_i}$ to $P_e$. Then map, for each $e\in E(F)$, the vertices $y_{e_{i+\ell}},y_{e_{-(i+\ell)}}$ to $X_{\text{final}}$ in an arbitrary injective manner. Also, map each $y_{e_{i+2\ell}},y_{e_{-(i+2\ell)}}$ to $U$ in an arbitrary injective manner. More generally, map the vertices $y_{e_{i+a\ell}},y_{e_{-(i+a\ell)}}$ with $a\geq 1$ odd to $X_{\text{final}}$ in an arbitrary injective manner and map the vertices $y_{e_{i+a\ell}},y_{e_{-(i+a\ell)}}$ with $a\geq 2$ even to $U$ in an arbitrary injective manner. We then just need to find paths of length $\ell$ connecting $y_{e_{i+a\ell}}$ and $y_{e_{i+(a+1)\ell}}$ (and paths of length $\ell$ connecting $y_{e_{-(i+a\ell)}}$ and $y_{e_{-(i+(a+1)\ell)}}$) which are disjoint from each other and from the images of the already mapped vertices. Since $(x,y)$ is $(\ell,L)$-bad for every $x\in X_{\text{final}}$ and $y\in U\cup \bigcup_{e\in E(F)} R_e$, such paths exist by Lemma \ref{manydisjoint}, provided that $L$ is sufficiently large.
\end{proof2}

\begin{corollary} \label{admissiblenotgood}
	Let $G$ be an $H$-free $K$-almost-regular graph on $n$ vertices with minimum degree~$\d\geq L^{100^k|V(H)|}$. Then, provided that $L$ is sufficiently large compared to $|V(H)|$ and $K$, for any $2\leq \ell\leq k$, the number of $L$-admissible, but not $L$-good, paths of length $\ell$ is at most $n\frac{2(K\d)^\ell}{f(\ell-1,L)}$.
\end{corollary}

Now we are in a position to prove Lemma \ref{kpathsgood}.

\begin{proof}[\textnormal{\textbf{Proof of Lemma \ref{kpathsgood}}}]
	Suppose that the path $u_0u_1\dots u_k$ is not $L$-good. Take $0\leq i<j\leq k$ with $j-i$ minimal such that $u_iu_{i+1}\dots u_j$ is not $L$-good. Then $u_i\dots u_j$ is $L$-admissible. For any fixed $i,j$, by Corollary \ref{admissiblenotgood}, the number of such paths is at most $n\frac{2(K\d)^{j-i}}{f(j-i-1,L)}\cdot 2(K\d)^{k-(j-i)}=4K^k\frac{n\d^k}{f(j-i-1,L)}\leq 4K^k\frac{n\d^k}{L}$. Using that $i$ and $j$ can take at most $k+1$ values each, it follows that the number of not $L$-good paths of length $k$ is at most $(k+1)^24K^k\cdot \frac{n\d^k}{L}$.
\end{proof}



\section{Long paths} \label{seclongpaths}


In what follows, for a vertex $x\in V(G)$ and a nonnegative integer $i$, we write $\cP_i(x)$ for the set of directed paths of length $i$ starting at $x$. For an element $P\in \cP_i(x)$, we let $v(P)$ be the endpoint of the path $P$.

\begin{definition}
	Let $i,j$ be nonnegative integers with $i+j<2k$. Call a pair $(x,y)$ of distinct vertices \emph{$(i,j)$-rich} if the number of pairs $(P,Q)\in \cP_i(x)\times \cP_j(y)$ such that there are at least $(|V(H)|+2)(2k+1)+1$ pairwise internally vertex-disjoint paths of length $2k-i-j$ between $v(P)$ and $v(Q)$ is more than $(2(i+j)|V(H)|(2k+1)+2(i+1)j)(K\d)^{i+j-1}$. Otherwise (including when $x=y$) call it \emph{$(i,j)$-poor}.
\end{definition}

\begin{lemma} \label{richtopaths}
	Let $G$ be a graph with maximum degree at most $K\d$. Let $x,y\in V(G)$ and let $i,j$ be nonnegative integers with $i+j<2k$. If $(x,y)$ is $(i,j)$-rich, then there exist $|V(H)|$ pairwise internally vertex-disjoint paths of length $2k$ between $x$ and $y$.
\end{lemma}

\begin{proof2}
	Choose a maximal set of pairwise internally vertex-disjoint paths $R_1,\dots,R_{\alpha}$ between $x$ and $y$ and assume that $\alpha<|V(H)|$. Let $T$ be the set of the vertices appearing in at least one of these paths. Note that $|T|<|V(H)|(2k+1)$.
	
	\smallskip
	
	\noindent \emph{Claim.} If there is a pair $(P,Q)\in \cP_i(x)\times \cP_j(y)$ such that
	\begin{enumerate}[label=(\roman*)]
		\item $P$ is disjoint from $T\setminus \{x\}$
		\item $Q$ is disjoint from $T\setminus \{y\}$
		\item $P$ and $Q$ are vertex-disjoint and
		\item there are at least $(|V(H)|+2)(2k+1)+1$ pairwise internally vertex-disjoint paths of length $2k-i-j$ between $v(P)$ and $v(Q)$,
	\end{enumerate} then there is a path of length $2k$ between $x$ and $y$ which is internally vertex-disjoint from all of $R_1,\dots,R_{\alpha}$.

	\smallskip

	\noindent \emph{Proof of Claim.} Clearly, it suffices to find a path of length $2k-i-j$ between $v(P)$ and $v(Q)$ which is disjoint from the vertices of $R_1,\dots,R_{\alpha},P,Q$, except for $v(P)$ and $v(Q)$. But such a path exists since there are at most $(\alpha+2)\cdot (2k+1)\leq (|V(H)|+2)(2k+1)$ vertices in one of $R_1,\dots,R_{\alpha},P,Q$ and there are at least $(|V(H)|+2)(2k+1)+1$ pairwise internally vertex-disjoint paths of length $2k-i-j$ between $v(P)$ and $v(Q)$.
	
	\smallskip
	
	A path provided by the claim would contradict the maximality of $R_1,\dots,R_{\alpha}$, so it suffices to prove that there are paths $P,Q$ satisfying (i)-(iv) above.
	
	Since the maximum degree of $G$ is at most $K\d$, the number of paths of length $i-1$ in $G$ intersecting $T$ is at most $i|T|(K\d)^{i-1}$, so the number of $P\in \cP_i(x)$ which have a vertex in $T\setminus \{x\}$ is at most $2i|T|(K\d)^{i-1}$. Since $|\cP_j(y)|\leq (K\d)^j$, the number of pairs $(P,Q)\in \cP_i(x)\times \cP_j(y)$ failing condition (i) above is at most $2i|T|(K\d)^{i-1}(K\d)^j$. Similarly, the number of pairs failing (ii) is at most $2j|T|(K\d)^{j-1}(K\d)^i$. Finally, for every $P\in \cP_i(x)$, the number of paths of length $j-1$ which intersect $P$ is at most $(i+1)j(K\d)^{j-1}$, so the number of pairs $(P,Q)\in \cP_i(x)\times \cP_j(y)$ for which $P$ and $Q$ share a vertex other than $y$ is at most $(K\d)^{i}\cdot 2(i+1)j(K\d)^{j-1}$. So the number of pairs which fail at least one of (i),(ii),(iii) is at most $(2(i+j)|T|+2(i+1)j)(K\d)^{i+j-1}\leq (2(i+j)|V(H)|(2k+1)+2(i+1)j)(K\d)^{i+j-1}$. By the definition of $(i,j)$-richness of $(x,y)$ it follows that there is a pair $(P,Q)$ satisfying (i)-(iv).
\end{proof2}

\begin{definition}
	For a vertex $v\in V(G)$ and some $1\leq \ell\leq k$, define an auxiliary graph $\mathcal{G}_{\ell}(v)$ as follows. The vertices of $\mathcal{G}_{\ell}(v)$ are the $(k+1)$-tuples $(u_0,u_1,\dots,u_k)\in V(G)^{k+1}$ with $u_0=v$ such that $u_iu_{i+1}\in E(G)$ for all $i$. Vertices $(u_0,\dots,u_k)$ and $(u'_0,\dots,u'_k)$ are joined by an edge if $v,u_1,u_2,\dots,u_k,u'_1,\dots,u'_k$ are distinct and there exist $0\leq i,j\leq k-1$ such that the pair $(u_{\ell},u'_{\ell})$ is $(i,j)$-rich. Since the vertex set of $\mathcal{G}_{\ell}(v)$ does not depend on $\ell$, we may define $\mathcal{G}(v)$ to be the union $\bigcup_{1\leq \ell\leq k} \mathcal{G}_{\ell}(v)$.
\end{definition}

\begin{lemma} \label{auxkrfree}
	Let $G$ be a graph with maximum degree at most $K\d$ which does not contain $H$ as a subgraph. Let $t=|V(F)|$. Then for any $v\in V(G)$ and any $1\leq \ell\leq k$, the graph $\mathcal{G}_{\ell}(v)$ is $K_t$-free.
	
	Moreover, let $r=R_k(t)$ be the $k$-colour Ramsey number. Then $\mathcal{G}(v)$ is $K_r$-free.
\end{lemma}

\begin{proof2}
	Suppose that $\mathcal{G}_\ell(v)$ contains $K_t$ as a subgraph. Let the corresponding vertices be the vectors $u^1,\dots,u^t$. Let their respective $(\ell+1)$th coordinate be $u^1_\ell,\dots,u^t_\ell$. For every $a\neq b$, since $u^{a}u^{b}$ is an edge in $\mathcal{G}_\ell(v)$, it follows that $u^a_\ell$ and $u^b_\ell$ are distinct, and, by Lemma \ref{richtopaths}, there exist $|V(H)|$ pairwise internally vertex-disjoint paths of length $2k$ between them. It is not hard to see that this implies that there is a copy of $H$ in $G$ in which the vertices of $F$ are mapped to $u^1_\ell,\dots,u^t_\ell$. This is a contradiction, so $\mathcal{G}_\ell(v)$ is indeed $K_t$-free.
	
	Suppose there is a copy of $K_r$ in $\mathcal{G}(v)$. Then each edge in this $K_r$ can be coloured with one of the colours $1,2,\dots,k$ such that if an edge gets colour $i$, then it lies in $\mathcal{G}_i(v)$. By the definition of $r$, there exists a monochromatic $K_t$ in this $k$-edge-coloured $K_r$, which gives a $K_t$ in some $\mathcal{G}_\ell(v)$, contradicting the first paragraph. 
\end{proof2}

\medskip

The next lemma provides us a large set of walks of length $2k$ with both endpoints in $S$. Later, we will argue that most of them are $L$-good paths.

\begin{lemma} \label{manypoor}
	Let $r=R_k(t)$ denote the $k$-colour Ramsey number where $t=|V(F)|$. Let $G$ be an $H$-free $K$-almost-regular graph on $n$ vertices with minimum degree $\d$ and let $S\s V(G)$ such that $|S|\geq 2nr/\d^{k}$. Then there are at least $\frac{|S|^2\d^{2k}}{4r^2n}$ vectors $(u_{-k},\dots,u_k)\in V(G)^{2k+1}$ with the following properties
	\begin{enumerate}[label=(\roman*)]
		\item $u_{-k}\in S$, $u_k\in S$
		\item $u_\ell u_{\ell+1}\in E(G)$ for every $-k\leq \ell\leq k-1$
		\item $(u_{-\ell},u_\ell)$ is $(i,j)$-poor for every $1\leq\ell\leq k$ and every $0\leq i,j\leq k-1$.
	\end{enumerate}
\end{lemma}

\begin{proof2}
	Since the minimum degree of $G$ is $\d$, the number of $(k+1)$-tuples $(v_0,v_1,\dots,v_k)\in V(G)^{k+1}$ with $v_{k}\in S$ and $v_iv_{i+1}\in E(G)$ for every $0\leq i\leq k-1$ is at least $|S|\d^k$. Writing $\mathcal{T}(v_0)$ for the set of such vectors for a fixed $v_0$ and letting $g(v_0)=|\mathcal{T}(v_0)|$, we get that $\sum_{v_0\in V(G)} g(v_0)\geq |S|\d^k$. Note that $\sum_{v_0\in V(G): g(v_0)<r} g(v_0)\leq nr\leq \frac{|S|\d^{k}}{2}$, so
	\begin{equation}
		\sum_{v_0\in V(G): g(v_0)\geq r} g(v_0)\geq \frac{|S|\d^{k}}{2}. \label{eqngv}
	\end{equation}
	Note that $\mathcal{T}(v_0)\s V(\mathcal{G}(v_0))$. By Lemma \ref{auxkrfree}, the graph $\mathcal{G}(v_0)\lbrack \mathcal{T}(v_0)\rbrack$ is $K_r$-free. This graph has $g(v_0)$ vertices, so if $g(v_0)\geq r$, then the number of non-edges in $\mathcal{G}(v_0)\lbrack \mathcal{T}(v_0)\rbrack$ is at least $\frac{1}{{r\choose 2}}{g(v_0)\choose 2}\geq \frac{g(v_0)^2}{r^2}$. But if $v=(v_0,v_1,\dots,v_k)\in \mathcal{T}(v_0)$ and $v'=(v_0,v'_1,\dots,v'_k)\in \mathcal{T}(v_0)$ are such that $vv'$ is not an edge in $\mathcal{G}(v_0)$, then $(u_{-k},\dots,u_k)=(v'_k,v'_{k-1},\dots,v'_1,v_0,v_1,\dots,v_k)$ satisfies all three properties in the statement of the lemma. Therefore the number of such $(2k+1)$-tuples with $u_0=v_0$ is at least $\frac{g(v_0)^2}{r^2}$ provided that $g(v_0)\geq r$. By (\ref{eqngv}) and Jensen's inequality, we get $\sum_{v_0\in V(G): g(v_0)\geq r} \frac{g(v_0)^2}{r^2}\geq \frac{|S|^2\d^{2k}}{4r^2n}$, and the proof is complete.
\end{proof2}


\medskip

The following simple lemma shows that most walks of length $2k$ are paths.

\begin{lemma} \label{tuplesarepaths}
	Let $G$ be a graph on $n$ vertices with maximum degree at most $K\d$. Then the number of $(2k+1)$-tuples $(u_{-k},\dots,u_k)\in V(G)^{2k+1}$ such that $u_iu_{i+1}\in E(G)$ for every $i$ and $u_i=u_j$ for some $i\neq j$ is at most ${2k+1 \choose 2}K^{2k-1}\cdot n\d^{2k-1}$.
\end{lemma}

\begin{proof2}
	There are ${2k+1 \choose 2}$ ways to choose the pair $\{i,j\}$ and there are $n$ ways to choose $u_i=u_j$. Given any such choices, there are at most $(K\d)^{2k-1}$ ways to choose the vertices $u_b$ for $b\not\in \{i,j\}$ since any vertex in $G$ has degree at most $K\d$.
\end{proof2}

\medskip

Our strategy now is to take all the paths guaranteed by Lemmas \ref{manypoor} and \ref{tuplesarepaths} and discard those which contain a subpath of length $k$ which is not $L$-good. The next result shows that doing this we discard only a small proportion of the paths.

\begin{lemma} \label{ksubpathsgood}
	Let $G$ be an $H$-free $K$-almost-regular graph on $n$ vertices with minimum degree $\d\geq L^{100^k|V(H)|}$. Then, provided that $L$ is sufficiently large compared to $|V(H)|$ and $K$, the number of paths $u_{-k}u_{-k+1}\dots u_k$ of length $2k$ in $G$ with the property that there is some $-k\leq j\leq 0$ for which the path $u_ju_{j+1}\dots u_{j+k}$ is not $L$-good is $O(\frac{n\d^{2k}}{L})$.
\end{lemma}

\begin{proof2}
	By Lemma \ref{kpathsgood}, there are $O(\frac{n\d^k}{L})$ paths $u_ju_{j+1}\dots u_{j+k}$ which are not $L$-good, and since the maximum degree of $G$ is at most $K\d$, there are at most $2(K\d)^k$ ways to extend such a path to a path $u_{-k}u_{-k+1}\dots u_k$ of length $2k$. The result follows after summing these terms for all $-k\leq j\leq 0$.
\end{proof2}

\medskip

The next lemma is the first step to relate the notion of $L$-goodness with the notion of $(i,j)$-richness.

\begin{lemma} \label{2kbadkgood}
	Suppose that $u_{-k}u_{-k+1}\dots u_k$ is a path in $G$ which is not $L$-good but each of its subpaths of length $k$ is $L$-good. Then, provided that $L$ is sufficiently large compared to $|V(H)|$, there exist $1\leq \alpha,\beta\leq k$ with $\alpha+\beta> k$ such that there exist $(|V(H)|+2)(2k+1)+1$ pairwise internally vertex-disjoint paths of length $\alpha+\beta$ between $u_{-\alpha}$ and $u_{\beta}$.
\end{lemma}

\begin{proof2}
	Choose $-k\leq i<j\leq k$ with $j-i$ minimal such that $u_iu_{i+1}\dots u_j$ is not $L$-good. By the minimality of $j-i$, every proper subpath of $u_iu_{i+1}\dots u_j$ is $L$-good, so $u_iu_{i+1}\dots u_j$ is $L$-admissible. By Lemma \ref{manydisjoint}, there exist $(|V(H)|+2)(2k+1)+1$ pairwise internally vertex-disjoint paths of length $j-i$ between $u_i$ and $u_j$.
	
	By the assumption that every subpath of $u_{-k}u_{-k+1}\dots u_k$ of length $k$ is good, we have $j-i>k$, so $i<0$ and $j>0$. Thus, the choices $\alpha=-i$ and $\beta=j$ satisfy the conditions described in the lemma.
\end{proof2}

\medskip

The next result is the final ingredient to the proof of Lemma \ref{good2k}.

\begin{lemma} \label{poorgivefew}
	Let $G$ be a graph on $n$ vertices with maximum degree at most $K\d$. Then there are $O(n\d^{2k-1})$ paths $u_{-k}u_{-k+1}\dots u_k$ in $G$ with the following two properties
	\begin{enumerate}[label=(\roman*)]
		\item $(u_{-\ell},u_{\ell})$ is $(i,j)$-poor for every $1\leq \ell\leq k$ and every $0\leq i,j\leq k-1$ and
		\item there exist $1\leq \alpha,\beta\leq k$ with $\alpha+\beta> k$ such that there exist $(|V(H)|+2)(2k+1)+1$ pairwise internally vertex-disjoint paths of length $\alpha+\beta$ between $u_{-\alpha}$ and $u_{\beta}$.
	\end{enumerate}
\end{lemma}

\begin{proof2}
	Fix a pair $(\alpha,\beta)$ with $1\leq \alpha,\beta\leq k$ and $\alpha+\beta> k$. It suffices to prove that the number of paths satisfying (i) and (ii) for \emph{this} pair $(\alpha,\beta)$ is $O(n\d^{2k-1})$.
	
	Let $\ell=\alpha+\beta-k$. Note that $1\leq \ell\leq k$. Also, let $i=\alpha-\ell=k-\beta$ and $j=\beta-\ell=k-\alpha$. Observe that $0\leq i,j\leq k-1$.
	
	Suppose that $u_{-\ell}u_{-\ell+1}\dots u_\ell$ is a path such that $(u_{-\ell},u_\ell)$ is $(i,j)$-poor. By the definition of $(i,j)$-poorness, the number of pairs of paths $(u_{-\ell}u_{-\ell-1}\dots u_{-\alpha},u_{\ell}u_{\ell+1}\dots u_{\beta})$ such that there exist $(|V(H)|+2)(2k+1)+1$ pairwise internally vertex-disjoint paths of length $\alpha+\beta=2k-i-j$ between $u_{-\alpha}$ and $u_{\beta}$ is $O(\d^{i+j-1})$. Thus, the number of ways to extend $u_{-\ell}u_{-\ell+1}\dots u_\ell$ to a path  $u_{-k}u_{-k+1}\dots u_k$ possessing property (ii) with our fixed choice of $\alpha$ and $\beta$ is $O(\d^{i+j-1}\cdot (K\d)^{k-\alpha+k-\beta})=O(\d^{2k-2\ell-1})$, where the first factor bounds the number of possible ways to extend to $u_{-\alpha}u_{-\alpha+1}\dots u_{\beta}$, and the second factor bounds the number of possible ways to extend that to $u_{-k}u_{-k+1}\dots u_k$. The number of possible choices for $u_{-\ell}u_{-\ell+1}\dots u_{\ell}$ is $O(n\d^{2\ell})$, so the result follows.
\end{proof2}

\medskip

We are now in a position to complete the proof of Lemma \ref{good2k}.

\begin{proof}[\textnormal{\textbf{Proof of Lemma \ref{good2k}}}]
	The condition $|S|=\omega(\frac{n}{\d^{1/2}})$ implies that $n\d^{2k-1}=o(\frac{|S|^2\d^{2k}}{n})$, so by Lemmas \ref{manypoor} and \ref{tuplesarepaths}, there are $\Omega(\frac{|S|^2\d^{2k}}{n})$ paths $u_{-k}u_{-k+1}\dots u_k$ with both endpoints in $S$ such that $(u_{-\ell},u_\ell)$ is $(i,j)$-poor for every $1\leq \ell\leq k$ and every $0\leq i,j\leq k-1$. Discard all those paths among these in which there is a subpath of length $k$ which is not $L$-good. By Lemma \ref{ksubpathsgood}, we discarded $O(\frac{n\d^{2k}}{L})$ paths, which is $o(\frac{|S|^2\d^{2k}}{n})$, by the condition $|S|=\omega(\frac{n}{L^{1/2}})$. Of the remaining paths, discard all those for which there exist $1\leq \alpha,\beta\leq k$ with $\alpha+\beta> k$ such that there exist $(|V(H)|+2)(2k+1)+1$ pairwise internally vertex-disjoint paths of length $\alpha+\beta$ between $u_{-\alpha}$ and $u_{\beta}$. By Lemma \ref{poorgivefew}, there are $O(n\d^{2k-1})$ such paths, which is again $o(\frac{|S|^2\d^{2k}}{n})$. Hence, we are left with $\Omega(\frac{|S|^2\d^{2k}}{n})$ paths.
	
	We claim that each such path is $L$-good. Suppose otherwise, and take a path $u_{-k}u_{-k+1}\dots u_k$ which is not $L$-good. Since each of its subpaths of length $k$ is $L$-good, by Lemma \ref{2kbadkgood} there exist $1\leq \alpha,\beta\leq k$ with $\alpha+\beta> k$ such that there exist $(|V(H)|+2)(2k+1)+1$ pairwise internally vertex-disjoint paths of length $\alpha+\beta$ between $u_{-\alpha}$ and $u_{\beta}$. But we discarded these paths, which is a contradiction, and the proof is complete.
\end{proof}

\bibliographystyle{abbrv}
\bibliography{biblongersub}


\end{document}